\documentclass[matprg]{svjour}
\usepackage{amsmath,graphicx,psfrag}
\newcommand{\reals}{\mathbf{R}}

\newcommand{\Expect}{\mathbf{E}}
\newcommand{\ie}{{\em i.\/e.}}

\renewenvironment{proof}{\textbf{Proof.}}{\qed \bigskip}
\begin{document}

\author{Alexandre d'Aspremont \and Laurent El Ghaoui}
\title{Static Arbitrage Bounds on Basket Option Prices}
\institute{Alexandre d'Aspremont \at ORFE, Princeton University, Princeton NJ 08544, USA.\\
\email{alexandre.daspremont@m4x.org} \and Laurent El Ghaoui \at
Department of Electrical Engineering and Computer Sciences, Cory
Hall, University of California, Berkeley, CA, 94720, USA.
\email{elghaoui@eecs.berkeley.edu}} \mail{Alexandre d'Aspremont}

\date{Received: date / Revised version: date}
\subclass{44A12, 44A60, 90C05, 90C34, 91B28}
\combirunning{d'Aspremont \& El Ghaoui: Static Arbitrage Bounds on Basket Option Prices}
\maketitle

\begin{abstract}
We consider the problem of computing upper and lower bounds on the
price of an European basket call option, given prices on other
similar options. Although this problem is hard to solve exactly in
the general case, we show that in some instances the upper and
lower bounds can be computed via simple closed-form expressions,
or linear programs. We also introduce an efficient linear
programming relaxation of the general problem based on an integral
transform interpretation of the call price function. We show that
this relaxation is tight in some of the special cases examined
before.
\end{abstract}

\keywords{Arbitrage, Linear Programming, {R}adon Transform, Basket
Options, Moment Problems.}

\begin{acknowledgement}
The authors are very grateful to two anonymous referees for
numerous comments and corrections. We also benefited from
discussions with Marco Avellaneda, Dimitris Bertsimas, Stephen
Boyd, Rama Cont, Darrell Duffie, Hans F\"ollmer, Nicole El Karoui,
Noureddine El Karoui, Francesco Rapisarda and seminar participants
at the groupe de recherche FIQAM meeting at the Ecole
Polytechnique, the INFORMS 2003 annual meeting in Atlanta and the
IMA workshop on Risk Management and Model Specification Issues in
Finance. The work of the second author is supported in part by a
National Science Foundation Career Award NSF-ECS-9983874. This
paper first appeared as ERL Technical report UCB/ERL M02/33 in
Nov. 2002 at U.C. Berkeley and on ArXiv as math.OC/0302243 in Feb.
2003.
\end{acknowledgement}

\subsection*{Notation}
\label{ss:notation}

For two $n$-vectors $x, y$, $x \geq y$ (resp.\ $x<y$) means $x_i
\geq y_i$ (resp.\ $x_i<y_i$), $i=1, \ldots, n$; $x_+$ denotes the
positive part of $x$, which is the vector with components $\max
(x_i, 0)$. $e$ is the $n$-vector with all components equal to one,
and $e_i$ is the $i$-th unit vector of $\mathbf{R}^n$. The set
$\mathbf{R}_+^n$ denotes the set of $n$-vectors with non-negative
components, and $\mathbf{R}_{++}^n$ its interior. The cone of
nonnegative measures with support included in $\mathbf{R}_+^n$ is
denoted by $\mathcal{K}$. For $w \in \mathbf{R}^m$, $K \in
\mathbf{R}$ and $g \in \mathbf{R}^{m+1}$, the notation $\langle g,
(w,K) \rangle$ denotes the scalar product $\tilde{g}^Tw +
g_{m+1}K$, where $\tilde{g}$ contains the first $m$ elements of
$g$.

\section{Introduction}
\label{s:intro}

\subsection{Problem setup}
\label{ss:problem-setup}

Let $p \in \mathbf{R}_+^m$, $K_0 \in \mathbf{R}_+$, $w_0 \in
\mathbf{R}_{++}^n$ and $K_i \in \mathbf{R}_+$, $w_i \in
\mathbf{R}^n_+$, for $i=1, \ldots, m$. We consider the problem of
computing upper (resp. lower) bounds on the price of a European
basket call option with maturity $T$, strike $K_0$ and weight
vector $w_0$:
\begin{equation}
\begin{array}{ll}
\mbox{max./min.} & {\mathbf{E}}_{\pi }(w_0^{T}x-K_0)_+ \\
\mbox{subject to} & {\mathbf{E}}_{\pi }(w_i^{T}x-K_i)_+=p_i,\quad
i=1, \ldots, m, \label{eq:constraints}
\end{array}
\end{equation}
over all probability distributions $\pi \in \mathcal{K}$ on the
asset price vector $x$, consistent with a given set of observed
prices $p_i$ of options on other baskets. Note that we implicitly
assume that all the options have the same maturity, and that,
without loss of generality, the risk-free interest rate is zero
(we compare prices in the forward market).

We seek non-parametric bounds, that is, we do not assume any
specific model for the underlying asset prices; our sole
assumption is the absence of a static or ``buy-and-hold''
arbitrage today (i.e. an arbitrage that only requires trading
today and at the option maturity). The primary objective of these
bounds is not to detect and exploit arbitrage opportunities in the
basket option market, illiquidity issues are likely to make these
opportunities hard to exploit. However, the data on basket prices
(index options in equity markets or swaptions in fixed income) is
\textit{very sparse} and traders often rely on intuitive guesses
to extrapolate the remaining points, using these prices to
calibrate models and evaluate more complex derivatives. Our
results aim to provide an efficient method to check the validity
of these extrapolated prices where they are the most likely to
create static arbitrage opportunities, i.e. very far in or out of
the money.

From a financial point of view, our approach can be seen as a
one-period, non-parametric computation of the upper and lower
hedging prices defined in \cite{ElKa91} and \cite{Karou95} (see
also \cite {Kara98}). The necessary conditions we detail in
section \ref{s-general} have been extensively used in the
unidimensional case to infer information on the state-price
density given option prices (see \cite{Bree78} or \cite{Laur00}
among others), we study here a multidimensional generalization.

From an optimization point of view, problems such as the one above
have received a significant amount of attention in various forms.
First, we can think of problem (\ref{eq:constraints}) as a
\textit{linear semi-infinite program}, \emph{i.\/e.}\ a linear
program with a finite number of linear constraints on an infinite
dimensional variable. We use this interpretation and the related
duality results to compute closed-form solutions in some
particular cases. Secondly, we can see (\ref{eq:constraints}) as a
\textit{generalized moment problem}. This approach was
successfully used in dimension one by \cite{Bert00}, who solve the
one dimensional problem completely and show that the
multidimensional extension is NP-Hard. However, their relaxation
algorithm requires the solution of a number of linear programs
that is potentially exponential in $n$, the number of assets. This
makes the method prohibitive for large-scale problems. Finally, as
in \cite{henk90}, one can think of (\ref {eq:constraints}) as an
\textit{integral transform inversion problem}. This is the
approach we adopt to design an efficient relaxation in the general
case.

The contribution of this work is twofold. First, exploiting the
necessary convexity of arbitrage free call prices
\[
C(w,K)={\mathbf{E}}_{\pi}(w^{T}x-K)_+,
\]
where $C(w,K)$ is the price of a basket call option with weights
$w$ and strike $K$, we detail a relaxation technique providing
upper (resp. lower) bounds on the solution to
(\ref{eq:constraints}). The resulting infinite dimensional linear
program on the call price function can be solved exactly. Compared
to the decomposition method proposed by \cite{Bert00}, this
relaxation technique has the advantage of being polynomial-time in
the number of assets and constraints.

Secondly, in some particular cases, we provide exact solutions to
(\ref{eq:constraints}) that have a polynomial complexity in the
number of assets and constraints. We also obtain expressions for
the corresponding pricing measures, and use them later on to prove
tightness of the linear programming relaxation in the general
case. \cite{Laur03} and \cite{Laur04} derived independently an
equivalent upper bound in the same particular cases and exact
lower bounds in dimension 2.

Finally, using recent results by \cite{Hobs04} who compute an
optimal solution to (\ref{eq:constraints}) in the particular case
where one seeks an upper bound on the price of a basket option
with positive weights given only single asset option prices:
\begin{equation}
\begin{array}{ll}
\mbox{maximize} & {\mathbf{E}}_{\pi }(w_0^{T}x-K_0)_+ \\
\mbox{subject to} & {\mathbf{E}}_{\pi }(x_i-K^j_i)_+=p^j_i,\quad
i=1, \ldots, m,~j=1,\ldots,J^i, \label{eq:constraints-eq}
\end{array}
\end{equation}
we show that our linear programming relaxation is tight.

Our paper is organized as follows. We begin in section
\ref{ss:dual-general} by a brief reminder of the fundamental
duality between pricing measures and positive portfolios. In
section~\ref{s-general} we detail a relaxation for problem
(\ref{eq:constraints}) using shape constraints on the call price
$C(w,K)$ as a function of the weight vector $w$ and the strike
price $K$. Using linear programming duality arguments, we obtain
in section~\ref{s-indiv} closed-form formulas (or simple linear
programs) giving upper and lower bounds in some particular cases.
In section~\ref{s-general-tightness}, using \cite{Hobs04} we show
that the linear programming relaxation derived in section
\ref{s-general} is tight in the particular case
(\ref{eq:constraints-eq}) above. Finally,
section~\ref{s-numerical} provides some numerical examples.

\subsection{Semi-infinite programming duality}
\label{ss:dual-general} We begin by detailing a key duality result
linking the existence of a pricing measure (or state price density
in \cite{Duff96}) and the absence of an arbitrage portfolio. In
the general case, we can write the upper bound problem as a
semi-infinite program:
\begin{equation}
p^{\sup }:=\sup_{\pi \in \mathcal{K}}\int_{\mathbf{R}_+^n}\psi
(x)\pi (x)dx \mbox{ subject to }\int_{\mathbf{R}_+^n}\phi (x)\pi
(x)dx=p, \;\;\int_{\mathbf{R}_+^n}\pi (x)dx=1,  \label{eq:main-pb}
\end{equation}
where
\[
\psi (x):=(w_0^{T}x-K_0)_+, \;\;\phi _i(x):=(w_i^{T}x-K_i)_+,
\;\;i=1, \ldots, m.
\]
We define the Lagrangian (on $\mathcal{K}\times \mathbf{R}^{m+1}$):
\[
L(\pi, \lambda, \lambda _0)=\int_{\mathbf{R}_+^n}\psi (x)\pi (x)dx+\lambda
^T \left( p-\int_{\mathbf{R}_+^n}\phi (x)\pi (x)dx \right) +\lambda
_0(1-\int_{\mathbf{R}_+^n}\pi (x)dx),
\]
and, as in \cite{Hett93}, we can write the dual of
(\ref{eq:main-pb}) as:
\begin{equation}
\begin{array}{rcl}
d^{\mathrm{sup}} & := & \displaystyle \inf_{\lambda _0, \lambda
}:\ \lambda ^{T}p+\lambda _0~:~\lambda ^T \phi (x)+\lambda _0 \geq
\psi (x)\mbox{
for every }x \in \mathbf{R}_+^n \\
& = & \displaystyle \inf_{\lambda }:\ \sup_{x \geq 0}:\lambda ^{T}p+\psi
(x)-\lambda ^T \phi (x).
\end{array}
\label{eq:main-pb-dual}
\end{equation}
Both primal and dual problems have very intuitive financial
interpretations. The primal problem looks for a pricing measure
that maximizes the target option price while satisfying the
pricing constraints imposed by the current market conditions. The
dual problem looks for the least expensive portfolio of options
and cash, $\lambda ^T \phi (x)+\lambda _0$, that dominates the
option payoff $\psi (x)$. Of course, the dual problem above yields
an upper bound on the upper bound.

Similarly, the computation of the lower bound involves
\begin{equation}
p^{\inf }:=\inf_{\pi \in \mathcal{K}}\int_{\mathbf{R}_+^n}\psi
(x)\pi (x)dx\mbox{ subject to }\int_{\mathbf{R}_+^n}\phi (x)\pi
(x)dx=p, \;\;\int_{\mathbf{R}_+^n}\pi (x)dx=1,
\label{eq:main-pb-inf}
\end{equation}
whose dual is
\begin{equation}
\begin{array}{rcl}
d^{\mathrm{inf}} & := & \displaystyle \sup_{\lambda _0, \lambda
}:\lambda ^{T}p+\lambda _0~:~\lambda ^T \phi (x)+\lambda _0 \leq
\psi (x)\mbox{
for every }x \in \mathbf{R}_+^n \\
& = & \displaystyle \sup_{\lambda }:\inf_{x \geq 0}:\lambda ^{T}p+\psi
(x)-\lambda ^T \phi (x).
\end{array}
\label{eq:main-pb-dual-inf}
\end{equation}
Here, the dual problem provides a lower bound on the lower bound.

General results on semi-infinite linear programs establish the equivalence
between the primal and dual formulations. We cite here a sufficient
constraint qualification condition for perfect duality from \cite{Hett93},
which makes an assumption about the support of optimal distributions. (We
focus now on the lower bound; a similar result holds for the upper bound
problem.)

\begin{proposition}
Assume that in problem (\ref{eq:main-pb-dual-inf}), the support of
the asset price distribution can be restricted to a given compact
set $B \subset \mathbf{R}_+^n$. Assume further that there exists a
pair $(\lambda _0, \lambda )\in \mathbf{R}^{n+1}$ such that:
\[
\lambda ^T \phi (x)+\lambda _0<\psi (x)\mbox{ for all } x \in B .
\]
Then if $d^{\mathrm{inf}}$ is finite, perfect duality holds,
namely $d^{\mathrm{inf}}=$ $p^{\inf }$.
\end{proposition}
\begin{proof}
See \cite{Hett93}.
\end{proof}

This constraint qualification condition trivially holds when $\phi
(x)$ and $\psi (x)$ are Call option payoffs hence we have
$d^{\mathrm{inf}}=$ $p^{\inf} $, provided that the support of
distributions feasible for our problem can be restricted to some
compact $B \subset \mathbf{R}_+^n$. However, this is often not the
case for the bounds detailed in section~\ref{s-indiv} and we will
prove perfect duality directly whenever possible.

\section{Relaxation for the general case using an integral transform}
\label{s-general}
\subsection{The Radon transform} \label{ss:integral-transform}
Let us come back to problem (\ref{eq:constraints}), for $p \in
\mathbf{R}_+^m$, $K \in \mathbf{R}_+^m$, $w_0 \in
\mathbf{R}_{++}^n$, $w_i \in \mathbf{R}^n_{+} $, $i=1, ..., m$ and
$K_0 \geq 0$, we seek to compute upper and lower bounds on the
price of a European call basket option with strike $K_0$ and
weight vector $w_0$:
\[
{\mathbf{E}}_{\pi }(w_0^{T}x-K_0)_+,
\]
with respect to all probability distributions $\pi \in
\mathcal{K}$ on the asset price vector $x$, consistent with a
given set of $m$ observed prices $p_i $ of options on other
baskets and forward prices $q_i=\mathbf{E}_{\pi}x_j$, that is,
given
\[
{\mathbf{E}}_{\pi}(w_i^Tx-K_i)_+=p_i, \;\;i=1, \ldots, m ~\mbox{
and }~\mathbf{E}_{\pi}x_j=q_j, \;\;j=1, \ldots, n.
\]
Extending to basket options the results of \cite{Bree78}, we
write, for some $\pi \in \mathcal{K}$:
\begin{eqnarray*}
C(w, K) &=&{\mathbf{E}}_{\pi }(w^{T}x-K)_+ \\
&=&\int_{\mathbf{R}_+^n}(w^{T}x-K)_{+}d \pi (x),
\end{eqnarray*}
we can think of $C(w, K)$ as a particular integral transform of
the measure $\pi $ and we can try to compute its inverse. If we
assume that the measure $\pi $ is absolutely continuous with
respect to the Lebesgue measure with density $\pi (x)$, then for
almost all $K$ we have:
\[
\frac{\partial ^{2}C(w, K)}{\partial
K^2}=\int_{\mathbf{R}_+^n}\delta (w^{T}x-K)\pi (x)dx,
\]
where $\delta (x)$ is the Dirac Delta function. This means that
${\partial ^{2}C(w, K)}/{\partial K^2}$ is the Radon transform of
the measure $\pi$ (see \cite{Helga99} for example). In this
setting, the general pricing problem above can then be rewritten
as the following infinite dimensional problem:
\[
\begin{array}{ll}
\mbox{min./max.} & f(w_0, K_0) \\
\mbox{subject to} & f(w_j, K_j)=p_j, \quad j=1, \ldots, m \\
&f(e_i, 0)=q_i, \quad i=1,\ldots,n\\
& f(w, K)\in \mathcal{R}_C,
\end{array}
\]
where $e_i$ is the Euclidean basis in $\reals^n$ and
$\mathcal{R}_C$ is the range of the (linear) integral transform
\[
\begin{array}{ll}
C: & \mathcal{K \rightarrow R}_C \\
& \pi \mapsto C(w, K)=\int_{\mathbf{R}_+^n}(w^{T}x-K)_{+}d \pi
(x).
\end{array}
\]
Thus, the problem of finding all possible arbitrage-free option
prices becomes equivalent to that of characterizing the range of
the Radon transform on the set of nonnegative measures
$\mathcal{K}$. This has been done by \cite{henk90} in the context
of production functions (which can be thought of as Put options).
As in \cite{henk90}, we denote by $C^{\infty
}\{0,\mathbf{R}_+^n\}$ the set of functions $f$ such that  for any
$k$ there is a polynomial $P_k(x)$ of degree $k$ such that:
\[
f(x)-P_k(x)=o(|x|^k),\quad \mbox{as } |x|\rightarrow 0,~
x\in\mathbf{R}_+^n.
\]
Using Call-Put parity, we can directly derive from \cite[theorem
3.2]{henk90} the following result:

\begin{proposition}
\label{prop:nec-suff-convex} A function $C(w, K), $ with $w \in
\mathbf{R}_+^n$ and $K>0$ belongs to $\mathcal{R}_C$, i.e. it can
be represented in the form
\[
C(w, K)=\int_{\mathbf{R}_+^n}(w^{T}x-K)_{+}d \pi (x),
\]
where $\pi $ is a nonnegative measure on a compact of
$\mathbf{R}_+^n$, if and only if the following conditions hold.

\begin{itemize}
\item  $C(w, K)$ is convex and homogenous of degree one;

\item  for every $w \in \mathbf{R}_{++}^n$, we have
\[
\lim_{K \rightarrow \infty }C(w, K)=0 \mbox{ and }\lim_{K
\rightarrow 0^+}\frac{\partial C(w, K)}{\partial K}=-1;
\]

\item if we write $D_{\xi}=\sum_i \xi_i\partial/\partial x_i$, the function
\[
F(w)=\int_0^{\infty }e^{-K}d \left( \frac{\partial C(w,
K)}{\partial K}\right)
\]
belongs to $C^{\infty }\{0,\mathbf{R}_+^n\}$ and for some
$\tilde{w}\in \mathbf{R}_+^n$ the inequalities:
\[
\left( -1 \right) ^{k+1}D_{\xi _1}...D_{\xi _k}F(\lambda
\tilde{w})\geq 0
\]
hold for all positive integers $k$ and $\lambda \in
\mathbf{R}_{++}$ and all $\xi _1, ..., \xi _k$ in
$\mathbf{R}_+^n$.
\end{itemize}
\end{proposition}
\begin{proof}
See \cite{henk90}.
\end{proof}
This result generalizes the necessary conditions for the absence
of arbitrage used by \cite{Bree78}, \cite{Laur00} or \cite{Bert00}
in dimension one.

\subsection{Linear programming relaxation}
\label{ss:lp-relaxation} The conditions above are not tractable in
the general case but we can formulate a relaxation of the original
program by replacing the last (moment) condition with weaker
monotonicity and linearity conditions. We then get an upper bound
on the upper bound (resp. a lower bound on the lower bound)
solution by computing:
\begin{equation}
\begin{array}{ll}
\sup /\inf & C(w_{0},K_{0}) \\
\mbox{subject to} & C(w,K)\mbox{ (jointly) convex in }\left(w,K\right) \\
& C(w,K)\mbox{ homogeneous of degree }1 \\
& -1\leq \partial C(w,K)/\partial K\leq 0\mbox{ and }C(w,K)\mbox{
nondecreasing in }w \\
& C(w_{i},0)=w_{i}^{T}q\mbox{, }i=1,...,m \\
& C(w_{i},K_{i})=p_{i}\mbox{, }i=1,...,m.
\end{array}
\label{eq:infinite-LP}
\end{equation}
where the variable is here $C(w,K)\in C\left(
\mathbf{R}^{n+1}\rightarrow \mathbf{R}_{+}\right)$. As we show
below, this infinite program can be reduced to a finite LP. If we
define $p_{m+i}=w_{i}^{T}q$ and $K_{m+i}=0$ for $i=1,...,m$ and
$p_{2m+1}=w_0^{T}q$ \ with $K_{2m+1}=0$, we can show the following
result:

\begin{proposition}
\label{prop:relax-general} If the following finite LP:
\begin{equation}
\begin{array}{ll}
\text{\textrm{max./min.}} & p_{0} \\
\text{subject to} & \left\langle
g_{i},(w_{j},K_{j})-(w_{i},K_{i})\right\rangle \leq
p_{j}-p_{i},\quad
i,j=0,...,2m+1 \\
& g_{i,j}\geq 0,~-1\leq g_{i,n+1}\leq 0,\quad i=0,...,2m+1,\quad
j=1,...,n
\\
& \left\langle g_{i},(w_{i},K_{i})\right\rangle =p_{i},\quad
i=0,...,2m+1,
\end{array}
\label{eq:finite-LP}
\end{equation}
in the variables $p_{0}\in \mathbf{R}_{+}$ and $g_{i}\in
\mathbf{R}^{n+1}$ for $i=0,...,2m+1$, is strictly feasible and its
optimal value is finite, the infinite program
(\ref{eq:infinite-LP}) and its discretization (\ref
{eq:finite-LP}) have the same optimal value. Furthermore, an
optimal point of (\ref{eq:infinite-LP}) can be constructed from
the solution to (\ref{eq:finite-LP}).
\end{proposition}
\begin{proof}
As in \cite{Boyd03}, we first notice that as a discretization of
the infinite program (\ref{eq:infinite-LP}), the finite LP will
compute a lower (or upper) bound on its optimal value. Let us now
show that this bound is in fact equal to the optimal value of
(\ref{eq:infinite-LP}). If we note $z^{\ast }=\left[ p_0^{\ast },
g_0^{\ast T}, \ldots, g_k^{\ast T}\right] ^T$ the optimal solution
to the LP problem above and if we define:
\[
C(w, K)=\max_{i=0, \ldots, 2m+1}\left \{ p_i^{\ast }+\left \langle
g_i^{\ast }, (w, K)-(w_i, K_i)\right \rangle \right \},
\]
where $p_i^{\ast}=p_i$ for $i=1,\ldots,2m+1$. $C(w, K)$ satisfies
\[
C(w_{i},K_{i})=p_i, \quad i=1, \ldots, 2m+1,
\]
and, by construction, $C(w_0,K_0)$ attains the lower bound $p_0$
computed in the finite LP. Also, $C(w, K)$ is convex as the
pointwise maximum of affine functions and is piecewise affine with
gradient $g_i$, which implies that it also satisfies the convexity
and monotonicity conditions in (\ref {eq:infinite-LP}), hence it
is a feasible point of the infinite dimensional problem. This
means that both problems have the same optimal value and $C(w,K)$
is an optimal solution to the infinite dimensional program in
(\ref {eq:infinite-LP}).
\end{proof}

\section{Exact upper and lower bounds using LP duality}
\label{s-indiv} In this section, we address the problem of
computing exact bounds in some particular cases using linear
programming duality. We first consider the case when the observed
data set corresponds to option and forward prices on each
individual asset since in practice, observations always include
the forward contract prices $\mathbf{E}_{\pi}x_i=q_i$, $i=1,
\ldots, n$ (forward contracts exist whenever options do). The
expressions derived in these simple particular cases will be used
in the section \ref{s-general-tightness} to show tightness of the
upper bound relaxation in the general case.

We first examine the problem of computing upper and lower bounds
on
\[
{\mathbf{E}}_{\pi }(w^{T}x-K_{0})_{+},
\]
given the $2n$ constraints
\begin{equation}
{\mathbf{E}}_{\pi }(x_{i}-K_{i})_{+}=p_{i},\;\;\mathbf{E}_{\pi
}x_{i}=q_{i},\;\;i=1,\ldots ,n,  \label{eq:simple-constraints}
\end{equation}
where $K_{0}>0$ and $w,K,p,q$ are given vectors of
$\mathbf{R}_{++}^{n}$. We will assume that $0 \leq p < q \leq
p+K$, which is a necessary and sufficient condition for the
problem above to be feasible. We show sufficiency by constructing
a discrete asset price distribution that matches these prices. Let
us define marginal distributions $\pi_i(x_i)$ such that
\begin{equation}
\label{eq:dist-feas} x_i= \left \{
\begin{array}{ll}
0 & \mbox{with probability } \pi_i(0)=1-\frac{q_i-p_i}{K_i}\\
\frac{q_iK_i}{q_i-p_i} & \mbox{with probability } \pi_i(\frac{q_iK_i}{q_i-p_i})=\frac{q_i-p_i}{K_i}\\
\end{array}
\right.
\end{equation}
and because $0 \leq p < q \leq p+K$, we know that
\[0 < \frac{q_i-p_i}{K_i}\leq 1. \]
Using Sklar's theorem, we can then construct an asset distribution
$\pi(x)$ with marginals $\pi_i(x_i)$, \ie~matching the market
prices of single asset options, hence the market is arbitrage
free. From the form of the constraints in
(\ref{eq:simple-constraints}), we also observe that the
constraints $0\leq p < q\leq p+K$ are necessary.

\subsection{Upper bound}
\subsubsection{One forward and one option price constraint per
asset}
\label{s-upper-forwards} Here, we apply the semi-infinite
duality result to the upper bound problem described in
(\ref{eq:simple-constraints}) to show the following result:
\begin{proposition}
\label{prop:upper-forward} Let $p,~q \in \mathbf{R}_+^n$, $K_0 \in
\mathbf{R}_+$, $w \in \mathbf{R}_{++}^n$ and $K_i \in
\mathbf{R}_+$ for $i=1, \ldots, n$, with $0\leq p_i < q_i \leq p_i
+ K_i$. An upper bound on the optimal value of the problem:
\[
\begin{array}{ll}
\mbox{maximize} & {\mathbf{E}}_{\pi }(w^{T}x-K_0)_+ \\
\mbox{subject to} & {\mathbf{E}}_{\pi }(x_i)=q_i\\
& {\mathbf{E}}_{\pi}(x_i-K_i)_+=p_i,\quad i=1, \ldots, n,
\end{array}
\]
is given by:
\begin{equation}
d^{\mathrm{sup}}=\displaystyle \max_{0 \leq j \leq n+1}\:
w^{T}p+\sum_{i}w_i\displaystyle \min (q_i-p_i, \beta
_{j}K_i)-\beta _{j}K_0,
\end{equation}
with $\beta _j:=(q_j-p_j)/K_j \in [0,1]$, $j=1, \ldots, n$ and the
convention $\beta _0=0$, $\beta _{n+1}=1$.
\end{proposition}
\begin{proof}
In view of the general result (\ref{eq:main-pb-dual}), the dual problem can
be expressed as
\begin{equation}
\begin{array}{rcl}
d^{\sup } & = & \displaystyle \inf_{\lambda +\mu \geq
w}\:\displaystyle\sup_{x \geq 0}\;\lambda ^{T}p+\mu
^{T}q+(w^{T}x-K_0)_+-\lambda ^T(x-K)_+-\mu ^{T}x,
\end{array}
\label{eq:dual-up-simple}
\end{equation}
where, without loss of generality, we have included the constraint $\lambda
+\mu \geq w$, in order to ensure that the inner supremum is finite. We
introduce a partition of $\mathbf{R}_+^n$ as follows. To a given subset $I$
of $\{1, \ldots, n \}$, we associate a subset $D_I$ of $\mathbf{R}_+^n $,
defined by
\[
D_I=\left \{ x~:~x_i>K_i, \;i \in I, \;\;0 \leq x_i \leq K_i, \;i
\in I^c \right \},
\]
where $I^c$ denotes the complement of $I$ in $\{1, \ldots, n \}$.
For $z \in \mathbf{R}^n$, let $z_I$ be the vector formed with the
elements $(z_i)_{i \in I}$, in the ascending order of indices in
$I$.

We have for
\[
\begin{array}{lll}
d^{\mathrm{sup}} & = & \displaystyle\inf_{\lambda +\mu \geq
w}:\displaystyle\max_{t\in \{0,1\},I\subseteq \{1,\ldots
,n\}}:\sup_{x\in D_{I}}:\lambda
^{T}p+\mu ^{T}q \\
&  & +t(w^{T}x-K_{0})-\lambda _{I}^{T}(x_{I}-K_{I})-\mu ^{T}x \\
& = & \displaystyle\inf_{\lambda +\mu \geq w}:\displaystyle\max_{t\in
\{0,1\}I\subseteq \{1,\ldots ,n\}}:\lambda ^{T}p+\mu ^{T}q+h(\lambda ,\mu
,I,t),
\end{array}
\]
where $h(\lambda, \mu, I, t)$ is given by
\[
\begin{array}{lll}
h(\lambda ,\mu ,I,t) & := & \displaystyle\sup_{x\in
D_{I}}:t(w^{T}x-K_{0})-\lambda _{I}^{T}(x_{I}-K_{I})-\mu ^{T}x \\
& = & \displaystyle\sup_{0\leq x_{I^c}\leq K_{I^c}}:(tw_{I^c}-\mu
_{I^c})^{T}x_{I^c}-tK_{0}+\lambda _{I}^{T}K_{I} \\
&  & +\displaystyle\sup_{x_{I}>K_{I}}:(tw_{I}-\mu _{I}-\lambda
_{I})^{T}x_{I}
\\
& = & \left\{
\begin{array}{ll}
(tw_{I^c}-\mu _{I^c})_{+}^{T}K_{I^c}-tK_{0}+(tw_{I}-\mu _{I})^{T}K_{I} & \mbox{if }\lambda _{I}+\mu _{I}\geq tw_{I}, \\
+\infty  & \mbox{otherwise.}
\end{array}
\right.
\end{array}
\]
We note that finiteness of $h(\lambda, \mu, I, t)$ is guaranteed
by $\lambda +\mu \geq w$ and $0 \leq t \leq 1$. When these
conditions hold, the maximum value of $h(\lambda, \mu, I, t)$ over
$I \subseteq \{1, \ldots, n \}$ is obtained when the complement
$I^c$ is the full set, that is, when $I$ is empty. We obtain
\[
\begin{array}{rcl}
\displaystyle \max_{I \subseteq \{1, \ldots, n \}} \: h(\lambda, \mu, I, t)
& = & (tw-\mu )_+^{T}K-tK_0.
\end{array}
\]
Optimizing over $t$, we obtain
\[
\displaystyle \max_{t \in \{0, 1 \}}\:\displaystyle \max_{I \subseteq \{1,
\ldots , n \}} \: h(\lambda, \mu, I, t)=\max \left( (-\mu )_+^{T}K, (w-\mu
)_+^{T}K-K_0 \right) .
\]
This results in the following expression for $d^{\mathrm{sup}}$:
\begin{eqnarray}
d^{\mathrm{sup}} &=&\displaystyle \inf_{\lambda +\mu \geq w}\:\lambda
^{T}p+\mu ^{T}q+\max \left( (-\mu )_+^{T}K, (w-\mu )_+^{T}K-K_0 \right)
\nonumber \\
&=&\displaystyle \inf_{\mu }\: w^{T}p+\mu ^T(q-p)+\max \left( (-\mu
)_+^{T}K, (w-\mu )_+^{T}K-K_0 \right),  \label{eq:dsup-simple-mu}
\end{eqnarray}
which admits the following linear programming representation:
\[
\begin{array}{rcl}
d^{\mathrm{sup}}=\displaystyle \inf_{\mu, t, v, z}:w^{T}p+\mu ^T(q-p)+t & \:
& t \geq v^{T}K, \;\;v \geq 0, \;\;v+\mu \geq 0 \\
&  & t \geq z^{T}K-K_0, \;\;z \geq 0, \;\;z+\mu \geq w.
\end{array}
\]
The problem is feasible, and is thus equivalent to its dual. After
some elimination of dual variables, the dual problem can be
rewritten as:
\[
d^{\mathrm{sup}}=\max_{y, \beta }\:w^{T}p+w^{T}y-\beta K_0~:~
\begin{array}[t]{l}
(1-\beta )K \geq q-p-y \geq 0 \\
\beta K \geq y \geq 0.
\end{array}
\]
We remark that the above problem is feasible if and only if $p < q
\leq p+K$. We thus recover the primal feasibility condition
mentioned before. This condition ensures that the dual bound
$d^{\mathrm{sup}}$ is finite. The above further reduces to the
one-dimensional problem:
\begin{equation}
d^{\mathrm{sup}}=\displaystyle \max_{0 \leq \beta \leq
1}:w^{T}p+\sum_{i}w_i\displaystyle \min (q_i-p_i, \beta K_i)-\beta
K_0.  \label{eq:dsup-beta}
\end{equation}
The above problem is the maximization of a piecewise linear
concave function of one variable, thus the maximum is attained at
one of the break points $\beta _j:=(q_j-p_j)/K_j \in \lbrack
0,1]$, $j=1, \ldots, n$, or for $\beta =0, 1$. This way, we can
obtain a closed-form expression for the upper bound, namely
\[
d^{\mathrm{sup}}=\displaystyle \max_{0 \leq j \leq n+1}\:
w^{T}p+\sum_{i}w_i\displaystyle \min (q_i-p_i, \beta
_{j}K_i)-\beta _{j}K_0,
\]
with the convention $\beta _0=0$, $\beta _{n+1}=1$, which is the
desired result.
\end{proof}
Remark that, when the price of forwards is not given , the upper
bound is readily obtained by setting the variable $\mu$, which is
the variable dual to the constraint $\mathbf{E}_{\pi }x=q$, to
zero in the expression (\ref{eq:dsup-simple-mu}). We get the
simple closed-form expression
\begin{equation}
d^{\mathrm{sup}}=w^{T}p+(w^{T}K-K_0)_+,
\label{eq:upper-no-forward-bnd}
\end{equation}
which can be obtained as a direct consequence of Jensen's
inequality applied to the function $x \rightarrow x_+$. We can
check that the above bound satisfies some basic properties: it is
convex in $w$ and concave in $p, q$. Also, when $w=e_i$ (the
$i$-th unit vector), and $K_0=K_i$, we obtain
$d^{\mathrm{sup}}=p_i$, while for $K_i=0$ (\ie~the options are in
fact forward contracts), we obtain $d^{\mathrm{sup}}=p_i(=q_i)$.

\subsubsection{Two option price constraints per
asset} \label{ss:two-options} Using \cite{Hobs04}, the result we
just obtained can directly be extended to the (slightly) more
general case where two option price constraints are given for each
asset (but no forward price is specified). We will use this
technical result to show tightness of the LP relaxation in section
\ref{s-general-tightness}.

We let $p^1,p^2 \in \mathbf{R}_+^n$, $K_0 \in \mathbf{R}_+$, $w
\in \mathbf{R}_{++}^n$ and $K^1,K^2 \in \mathbf{R}_+^n$ (with
$K^2>K^1$). We look for an upper bound on the optimal value of the
problem:
\begin{equation}
\begin{array}{ll} \label{eq:two-options}
\mbox{maximize} & {\mathbf{E}}_{\pi }(w^{T}x-K_0)_+ \\
\mbox{subject to} & {\mathbf{E}}_{\pi }(x_i)=q_i\\
& {\mathbf{E}}_{\pi}(x_i-K^j_i)_+=p_i^j,\quad i=1, \ldots,
n,~j=1,2.
\end{array}
\end{equation}
where $q\in\reals^n_+$ is here a \emph{variable}. From
(\ref{eq:infinite-LP}), to preclude arbitrage between options and
forwards, we must impose:
\[
\underline{q}_i:=p^1_i+K^1_i\frac{p^1_i-p^2_i}{K^2_i-K^1_i} \leq
q_i \leq p^1_i +K^1_i:=\overline{q}_i, \quad i=1,\ldots,n.
\]
As in \cite{Hobs04}, for each asset $x_i$, we can then form
$\bar{C}^{(i)}(K^j_i)$, the largest decreasing convex function
such that $\bar{C}^{(i)}(0)=p^1_i +K^1_i$ and
$\bar{C}^{(i)}(K^j_i)=p^j_i$ for $i=1,\ldots,n,~j=1,2$. Then every
maximal decreasing convex function $C^{(i)}(k)$ matching the
market prices can be written in terms of $q_i$ as:
\[
C^{(i)}(k)=\bar{C}^{(i)}(k)-(p^1_i+K^1_i-q_i)\frac{(K^1_i-k)_+}{K^1_i},
\quad i=1,\ldots,n.
\]
From \cite{Hobs04}, we know that for a given $q$ the upper bound
on (\ref{eq:two-options}) can be written:
\[
\inf_{\{\lambda\geq0,~\lambda^Te=1\}} \sum_{i=1}^n w_i
C^{(i)}(\frac{\lambda_i}{w_i}K_0)
\]
this means that an upper bound on the solution to
(\ref{eq:two-options}) for all possible values of $q$ can be found
by solving:
\[
\sup_{\{\underline{q}_i\leq q_i \leq \overline{q}_i\}}
\inf_{\{\lambda\geq0,~\lambda^Te=1\}} \sum_{i=1}^n w_i
\bar{C}^{(i)}\left(\frac{\lambda_i}{w_i}K_0\right)-(p^1_i+K^1_i-q_i)\frac{(K^1_i-(\frac{\lambda_i}{w_i}K_0))_+}{K^1_i}
\]
and because both $\bar{C}^{(i)}$ and
$(K^1_i-(\frac{\lambda_i}{w_i}K_0))_+$ are convex functions of
$\lambda$ (with $(p^1_i+K^1_i-q_i)\leq0$), this is also:
\[
\inf_{\{\lambda\geq0,~\lambda^Te=1\}} \sup_{\{\underline{q}_i\leq
q_i \leq \overline{q}_i\}} \sum_{i=1}^n w_i
\bar{C}^{(i)}\left(\frac{\lambda_i}{w_i}K_0\right)-(p^1_i+K^1_i-q_i)\frac{(K^1_i-(\frac{\lambda_i}{w_i}K_0))_+}{K^1_i}.
\]
The inner supremum is reached for $q_i=p^1_i+K_i$ and we get the
solution to (\ref{eq:two-options}) as:
\[
\inf_{\{\lambda\geq0,~\lambda^Te=1\}}\sum_{i=1}^n w_i
\bar{C}^{(i)}\left(\frac{\lambda_i}{w_i}K_0\right).
\]
If $q_i=p^1_i+K_i$, the corresponding measure in \cite{Hobs04}
places no weight on values of $x_i$ smaller than $K^1_i$, hence
the problem with two options:
\[
\begin{array}{ll}
\mbox{maximize} & {\mathbf{E}}_{\pi }(w^{T}x-K_0)_+ \\
\mbox{subject to} & {\mathbf{E}}_{\pi}(x_i-K^j_i)_+=p_i^j,\quad
i=1, \ldots, n,~j=1,2.
\end{array}
\]
is equivalent to the following problem (setting $x=K^1+y$):
\[
\begin{array}{ll}
\mbox{maximize} & {\mathbf{E}}_{\pi }(w^{T}y-(K_0-w^TK^1))_+ \\
\mbox{subject to} & {\mathbf{E}}_{\pi }(y_i)=p^1_i\\
& {\mathbf{E}}_{\pi}(y_i-(K^2_i-K^1_i))_+=p^2_i,\quad i=1, \ldots,
n.
\end{array}
\]
where one forward and one option price constraint are given per
asset.

\subsection{Perfect duality: upper bound}
\label{s-upper-forwards-proba} We first compute the optimal
probability measures corresponding to the upper bound result with
option and forward price constraints obtained in
section~\ref{s-upper-forwards}. We can recover an optimal
distribution, or a sequence of distributions which achieve the
bound in the limit. This provides a direct proof of the fact that
$p^{\mathrm{sup}}=d^{\mathrm{sup}}$, i.e. that the upper bound
computed in Proposition~\ref{prop:upper-forward} is tight.

Without loss of generality, we assume $e^{T}w=1$. In
(\ref{eq:dsup-beta}) we obtained:
\[
d^{\mathrm{sup}}=\displaystyle \sup_{0 \leq \beta \leq
1}:w^{T}p+\sum_{i}w_i\displaystyle \min (q_i-p_i, \beta K_i)-\beta
K_0,
\]
which can be rewritten (the $\min $ is taken elementwise):
\[
\sup_{0 \leq \beta \leq 1}\:w^T \min \{q-\beta K_{0}e, p+\beta
\left( K-K_0 \right) \},
\]
or again:
\[
\sup_{0 \leq \beta \leq 1}\inf_{t \in \lbrack 0,
1]^m}\:w^T((e-t)(q-\beta K_{0}e)+t(p+\beta \left( K_i-K_0 \right)
)).
\]
Using LP duality we know that this is also equal to (with
$e^{T}w=1$):
\[
\inf_{t \in \lbrack 0, 1]^m}\sup_{0 \leq \beta \leq 1}\:\beta
\left( w^{T}tK-K_0 \right) +w^T(e-t)q+w^{T}tp.
\]
We express the above as
\[
\inf_{t \in \lbrack 0, 1]^m}\:w^T(e-t)q+w^{T}tp+\left( w^{T}tK-K_0
\right) _+.
\]
This problem can be solved exactly as a finite linear program, and
we obtain $t^{\ast}$ such that:
\begin{equation}
d^{\mathrm{sup}}=w^T \left( (e-t^{\ast })q+t^{\ast }p \right)
+\left( w^{T}t^{\ast }K-K_0 \right) _+.  \label{eq:dsup-synthetic}
\end{equation}
We recognize here the expression of the upper bound on the price
of a basket, where we are only given the following option price
constraints:
\[
{\mathbf{E}}_{\pi }(x_i-\hat{K}_i)_+=\hat{p}_i, \;\;i=1, \ldots,
n,
\]
where $\hat{K}:=t^{\ast }K$ and $\hat{p}:=(e-t^{\ast })q+t^{\ast
}p$. In other words, as in \cite{Hobs04}, this upper bound is
equal to the bound we would obtain given only one (synthetic) call
price $\hat{p}_i$ per asset corresponding to a strike price
$\hat{K}_i$.

Suppose first that $(w^{T}K\leq K_0)$, in this case we have
$t^{\ast}=e$ and $d^{\mathrm{sup}}=w^Tp$, and as in \cite[p.15,
Case 3]{Hobs04} we can construct a sequence of distributions
converging to the optimal value. If however $(w^{T}K > K_0)$ and
$p>0$, then using \cite[p.14, Case 2]{Hobs04}, we know that the
bound is attained by a distribution with finite support. Finally,
if $(w^{T}K > K_0)$ and $p_i=0$ for some index $i$, \cite[p.12,
Case 1]{Hobs04} shows that a similar result holds.

Finally, using the result above, we know that this result is still
valid when we replace the option and forward price constraints
with two option price constraints.

\subsection{Lower bound}
\label{ss:lower-bound} We can obtain a similar result for the
lower bound problem. In this case however the solution is not in
closed form and involves solving a (polynomial size) linear
program.
\begin{proposition}
\label{prop:lower-forward} Let $p,~q \in \mathbf{R}_+^n$, $K_0 \in
\mathbf{R}_+$, $w \in \mathbf{R}_{++}^n$ and $K_i \in
\mathbf{R}_+$ for $i=1, \ldots, n$, with $0\leq p_i < q_i \leq p_i
+ K_i$. A lower bound on the optimal value of the problem:
\[
\begin{array}{ll}
\mbox{minimize} & {\mathbf{E}}_{\pi}(w^{T}x-K_0)_+ \\
\mbox{subject to} & {\mathbf{E}}_{\pi}(x_i)=q_i\\
& {\mathbf{E}}_{\pi}(x_i-K_i)_+=p_i,\quad i=1, \ldots, n,
\end{array}
\]
can be computed by solving the following linear program:
\begin{equation}
\begin{array}{lll}
d^{\mathrm{inf}}= & \displaystyle\sup_{\lambda ,\mu ,\alpha
_{0},\ldots
,\alpha _{n}} & \lambda ^{T}p+\mu ^{T}(q-K)+h~ \\
& \text{subject to} & \lambda +\mu \leq w \\
&  & h\leq \alpha _{0}(w^{T}K-K_{0})-(\alpha _{0}w-\mu
)_{+}^{T}K,\;0\leq
\alpha _{0}\leq 1 \\
&  & \forall i~:~h\leq \alpha _{i}(w^{T}K-K_{0})-\sum_{j\neq
i}(\alpha
_{i}w_{j}-\mu _{j})_{+}K_{j} \\
&  & (\lambda _{i}+\mu _{i})_{+}/w_{i}\leq \alpha _{i}\leq 1,
\end{array}
\end{equation}
which has $3n+2$ variables and $3n+4$ constraints.
\end{proposition}
\begin{proof}
In the lower bound case, the dual problem is
\[
d^{\mathrm{inf}}=\sup_{\lambda +\mu \leq w}\:\inf_{x \geq 0}\:\lambda
^{T}p+\mu ^{T}q+(w^{T}x-K_0)_+-\lambda ^T(x-K)_+-\mu ^{T}x,
\]
where we exploited the fact that the inner infimum is $-\infty $
unless $\lambda +\mu \leq w$.

Let us use the same notation as before. We have
\[
\begin{array}{lll}
d^{\mathrm{inf}} & = & \displaystyle\sup_{\lambda +\mu \leq
w}:\displaystyle\min_{I\subseteq \{1,\ldots
,n\}}:\displaystyle\inf_{x\in D_{I}}:\lambda
^{T}p+\mu ^{T}q \\
&  & +(w^{T}x-K_{0})_{+}-\lambda _{I}^{T}(x_{I}-K_{I})-\mu ^{T}x \\
& = & \displaystyle\sup_{\lambda +\mu \leq w}:\displaystyle\min_{I\subseteq
\{1,\ldots ,n\}}:\lambda ^{T}p+\mu ^{T}q+h(\lambda ,\mu ,I),
\end{array}
\]
where
\[
\begin{array}{rcl}
h(\lambda, \mu, I) & = & \displaystyle \inf_{x, y_0}\: y_0-\lambda
_I^T(x_I-K_I)-\mu ^{T}x~:~x \in D_I, \: y_0 \geq w^{T}x-K_0, \;y_0 \geq 0 .
\end{array}
\]
We have by linear programming duality
\[
h(\lambda, \mu, I) = \sup \:(\alpha w-\mu )^{T}K-\alpha
K_0-(\alpha w_{I^c}-\mu _{I^c})_+^{T}K_{I^c} ~:~
\begin{array}[t]{l}
\alpha w_I - \lambda _I-\mu _I \geq 0 \\
0 \leq \alpha \leq 1
\end{array}
\]
Thus
\[
\begin{array}{rcl}
d^{\mathrm{inf}} & = & \displaystyle \sup_{\lambda +\mu \leq w}\:\lambda
^{T}p+\mu ^T(q-K)+\displaystyle \:\displaystyle \min_{I \subseteq \{1,
\ldots , n \}}\:f(\lambda, \mu, I),
\end{array}
\]
where
\[
f(\lambda, \mu, I):=\sup_{\underline{\alpha }(\lambda, \mu, I)\leq
\alpha \leq 1}\:\alpha (w^{T}K-K_0)-(\alpha w_{I^c}-\mu
_{I^c})_+^{T}K_{I^c},
\]
and
\[
\underline{\alpha }(\lambda, \mu, I):=\max_{i \in
I}\:\displaystyle \frac{(\lambda _i+\mu _i)_+}{w_i},
\]
with the convention that $\underline{\alpha }(\lambda, \mu, I)=0$ when $I$
is empty.

Let $I$ be a non-empty subset of $\{1, \ldots, n \}$. Let $i \in \arg
\max_{i \in I}(\lambda _i+\mu _i)_+/w_i$. We observe that
\[
\underline{\alpha }(\lambda, \mu, I)=\underline{\alpha }(\lambda, \mu , \{i
\}),
\]
and
\[
f(\lambda, \mu, I)\geq f(\lambda, \mu, \{i \}),
\]
which dramatically reduces the complexity of the minimization subproblem:
instead of computing the minimum over all $2^n$ sets $I \subseteq \{1,
\ldots, n \}$ it is sufficient to pick $I$ in the set of \emph{singletons}
of $\{1, \ldots, n \}$, or $I=\emptyset $. Therefore, the problem reads as a
linear program
\begin{equation}
\begin{array}{lll}
d^{\mathrm{inf}}= & \displaystyle\sup_{\lambda ,\mu ,\alpha _{0},\ldots
,\alpha _{n}} & \lambda ^{T}p+\mu ^{T}(q-K)+h~ \\
& \text{subject to} & \lambda +\mu \leq w \\
&  & h\leq \alpha _{0}(w^{T}K-K_{0})-(\alpha _{0}w-\mu )_{+}^{T}K,\;0\leq
\alpha _{0}\leq 1 \\
&  & \forall i~:~h\leq \alpha _{i}(w^{T}K-K_{0})-\sum_{j\neq i}(\alpha
_{i}w_{j}-\mu _{j})_{+}K_{j} \\
&  & (\lambda _{i}+\mu _{i})_{+}/w_{i}\leq \alpha _{i}\leq 1,
\end{array}
\label{eq:LP-with-forwards}
\end{equation}
and can be solved efficiently, since it has $O(n)$ constraints and variables.
\end{proof}

\subsection{Perfect duality: lower bound without forwards}
\label{ss:lower-bound-without} Here we study a particular case of
the lower bound problem discussed in
section~$\ref{ss:lower-bound}$ where we are not given information
on forward prices. We can't prove perfect duality in the setting
of section~\ref{ss:lower-bound} but prove it below in a more
restrictive case. Without information on the forward prices, we
simply set the dual variable $\mu$ to zero in expression
(\ref{eq:LP-with-forwards}) to obtain:
\begin{equation}  \label{eq:lp-lower-bnd-simple}
\begin{array}{rcl}
d^{\inf} & = & \displaystyle \sup_{0 \leq \xi \leq e} \: w^Tp \xi
+ h ~:~ h \leq 0, \; h \leq \xi_i (w_iK_i - K_0), \; 1 \leq i \leq
n .
\end{array}
\end{equation}

We note that $d^{\mathrm{inf}}$ can now be expressed as the
solution of a non-linear, convex optimization problem:
\begin{equation}
d^{\inf }=\sup_{\xi } \: w^Tp \xi -\max_{1 \leq i \leq n}\;\xi
_i(K_0-w_{i}K_i)_+~:~0 \leq \xi \leq e,
\label{eq:primal-ncvx-lower-bnd}
\end{equation}
or from its dual:
\begin{equation}
d^{\inf }=\inf_{\nu } \: \sum_{i=1}^n \left( p_{i}w_i-\nu
_i(K_0-w_{i}K_i)_+\right) _+~:~\nu ^{T}e=1, \;\nu \geq 0 .
\label{eq:dual-lp-lower-bnd}
\end{equation}

We can reduce again this optimization problem to a line search
over a scalar parameter, by elimination of the variable $\xi $. We
obtain
\[
d^{\inf }=\sum_{\{i~:~K_{i}w_i\geq K_0\}}p_{i}w_i+\sup_{v \geq
0}\:\sum_{\{i~:~K_{i}w_i<K_0\}}p_{i}w_i\min (1, \displaystyle
\frac{v}{K_0-K_{i}w_i})-v.
\]
The minimization above can be further reduced to a closed-form
expression by noting that the piecewise-linear function (of $v$)
involved has break points at $\gamma _i=K_0-K_{i}w_i$ (for $i$
such that $\gamma _i>0$) and $0$. We obtain as before a lower
bound on the minimization problem
\[
\begin{array}{ll}
\mbox{minimize}   & {\mathbf{E}}_{\pi}(w^{T}x-K_0)_+ \\
\mbox{subject to} & {\mathbf{E}}_{\pi}(x_i-K_i)_+=p_i,\quad i=1,
\ldots, n,
\end{array}
\]
as:
\[
\begin{array}{lll}
d^{\inf } &=&\displaystyle \sum_{\{i~:~K_{i}w_{i}\geq
K_{0}\}}p_{i}w_{i}
\label{eq:closed-form-lower-forward} \\
&&+\displaystyle \max_{\{j~:~K_{j}w_{j}<K_{0}\}}\left(
\displaystyle \sum_{\{i~:~K_{i}w_{i}<K_{0}\}}p_{i}w_{i}\min
(1,\displaystyle\frac{K_{0}-K_{j}w_{j}}{K_{0}-K_{i}w_{i}})-K_{0}+w_{j}K_{j}\right)_{+}
\end{array}
\]

Now, the linear programming expression
(\ref{eq:dual-lp-lower-bnd}) allows us to recover an optimal asset
price distribution or a sequence of distributions that are optimal
in the limit, as follows. Let $\nu $ be an optimal vector for
problem (\ref{eq:dual-lp-lower-bnd}). Let $\mathcal{I}$ be the set
of indices $i$ such that $K_0>w_iK_i$. We note that $i \not \in
\mathcal{I}$ implies $\nu _i=0$. For simplicity we assume that
$\mathcal{I}=\{1, \ldots, m \}$, where $0 \leq m \leq n$ (the
choice $m=0$ corresponding to empty $\mathcal{I}$).

First we examine the case when $m=0$, that is, $\mathcal{I}$ is
empty. In other words, $\min_i w(i)K_i \geq K_0$, and therefore
$d^{\mathrm{inf}} = p^Tw$. For $\epsilon>0$, we can define an
asset price distribution $\pi_\epsilon$ such that
\begin{equation}
x= \left \{
\begin{array}{ll}
\epsilon^{-1}p+K & \mbox{with probability } \pi_\epsilon=\epsilon\\
0& \mbox{with probability } \pi_=1-\epsilon\\
\end{array}
\right.
\end{equation}
which satisfies the pricing constraints and:
\[
\Expect_{\pi_\epsilon}[(w_0^Tx-K_0)_+]=w^Tp+\epsilon(K-K_0),
\]
we recover the lower bound by taking the limit when $\epsilon$
goes to zero.

Next, we assume $m\geq 1$. Let $\alpha =(n-m)/m$. For $\epsilon $
such that
\[
\epsilon <\alpha ^{-1}\min_{1\leq i\leq m}\nu _{i}(\neq 0),
\]
we define the vector $\nu (\epsilon )$ by
\[
\nu_i(\epsilon) = \left \{
\begin{array}{ll}
\nu_i-\alpha \epsilon & \mbox{if } 1 \leq i \leq m, \\
\epsilon & \mbox{otherwise.}
\end{array}
\right.
\]
Since $\epsilon$ is small enough, the vector $\nu(\epsilon)$
satisfies the constraints of problem (\ref{eq:dual-lp-lower-bnd}).
We now define a distribution $\pi_{\epsilon}$ on the asset price
vector $x$ as follows.
\[
x = x^\epsilon(i) \mbox{ with probability } \nu_i(\epsilon),
\]
where
\[
x_j^\epsilon(i) = \left \{
\begin{array}{ll}
\displaystyle \frac{p_j}{\nu_j(\epsilon)} + K_j & \mbox{if } j=i, \\
0 & \mbox{otherwise.}
\end{array}
\right.
\]
Note that $x_j^\epsilon(i)$ is always well-defined, since
$\nu_j(\epsilon)
>0 $ for every $j$.

Let us check that the distribution $\pi_{\epsilon}$ of asset
prices satisfies the constraints in (\ref{eq:constraints}). For
every $j$, $1 \leq j \leq n$, we have
\[
\begin{array}{rcl}
\mathbf{E}_{\pi_{\epsilon}}(x_j-K_j)_+ & = & \sum_{i=1}^n
\nu_i(\epsilon)
(x_j^\epsilon(i)-K_j)_+ \\
& = & \nu_j(\epsilon) (x_j^\epsilon(j)-K_j)_+ \\
& = & p_j .
\end{array}
\]
We also check that with this choice of asset price distribution,
the objective in (\ref{eq:constraints}) attains the lower bound
$d^{\mathrm{inf}}$, when we let $\epsilon \rightarrow 0$. We have
\[
\begin{array}{rcl}
{\mathbf{E}}_{\pi_{\epsilon}}(w^{T}x-K_0)_+ & = & \sum_{i=1}^n \nu
_i(w^{T}x^{\epsilon }(i)-K_0)_+ \\
& = & \sum_{i=1}^n \nu _i(\epsilon
)(\sum_{j=1}^{n}w_{j}x_j^{\epsilon
}(i)-K_0)_+ \\
& = & \sum_{i=1}^n \nu _i(\epsilon )(w_{i}x_i^{\epsilon }(i)-K_0)_+ \\
& = & \sum_{i=1}^n \nu _i(\epsilon )(w_i(\frac{p_i}{\nu _i(\epsilon )}+K_i)-K_0)_+ \\
& = & \sum_{i=1}^n(w_{i}p_i-\nu _i(\epsilon )(K_0-w_{i}K_i))_+.
\end{array}
\]
Letting $\epsilon \rightarrow 0$, we obtain
\[
\begin{array}{rcl}
\displaystyle \mbox{\textrm{lim}}_{\epsilon \rightarrow
0}{\mathbf{E}}_{\pi_{\epsilon}}(w^{T}x-K_0)_+ & = &
\sum_{i=1}^m(w_{i}p_i-\nu
_i(K_0-w_{i}K_i))_++\sum_{i=m+1}^{n}w_{i}p_i \\
& = & \sum_{i=1}^n(w_{i}p_i-\nu _i(K_0-w_{i}K_i)_+)_+ \\
& = & d^{\mathrm{inf}},
\end{array}
\]
as claimed. This shows that $d^{\mathrm{inf}}=p^{\mathrm{inf}}$
and that the lower bound computed in
(\ref{eq:closed-form-lower-forward}) is tight in the absence of
constraints on forward prices.

\section{Tightness of the integral transform based LP relaxation}
\label{s-general-tightness} In this section, using results from
\cite{Hobs04} and section \ref{s-indiv} we show that the linear
programming relaxation derived in
Proposition~\ref{prop:relax-general} is tight, \ie~yields the
exact solution to problem (\ref{eq:constraints}) in the particular
case considered in \cite{Hobs04} where only market prices of
\emph{single asset} options are given for \emph{many strikes}.
\begin{proposition}
\label{prop:tightness} Let $p^j_i \in \mathbf{R}_+$, $K^j \in
\mathbf{R}_+$ for $j=1,\ldots,J^i$, $K_0 \in \mathbf{R}_+$ and
$w_0 \in \mathbf{R}^n_{++}$ for $i=1, \ldots, m$. The problems
\begin{equation}
\label{eq:equity-prob}
\begin{array}{ll}
\mbox{maximize} & {\mathbf{E}}_{\pi }(w_0^{T}x-K_0)_+ \\
\mbox{subject to} & {\mathbf{E}}_{\pi }(x_i-K^j_i)_+=p^j_i,\quad
i=1, \ldots, n,~j=1,\ldots,J^i,
\end{array}
\end{equation}
and
\begin{equation}
\label{eq:lp-tight}
\begin{array}{ll}
\text{\textrm{maximize}} & p_{0} \\
\text{subject to} & \left\langle
g^i_{k},(w^j_{l},K^j_{l})-(w^i_{k},K^i_{k})\right\rangle \leq
p^j_{l}-p^i_{k}\quad i,j=0,
\ldots, n,~k,l=1,\ldots,J^i\\
& g_{i,i}\geq 0,~-1\leq g_{i,n+1}\leq 0, \quad k=i,\ldots,n\\
& \left\langle g_{i},(w_{i},K_{i})\right\rangle =p_{i},\quad i=0,
\ldots, n,~j=1,\ldots,J^i\\
\end{array}
\end{equation}
where $w^j_i=e_i$ for $j=1,\ldots,J^i$ and $J^0=2$, have the same
optimal value.
\end{proposition}
\begin{proof} We first focus on the particular case where only
the forward price and one option price are given per asset,
\ie~$K^1_i=0$ and $J^i=2$ for $i=1,\ldots,n$. As in section
\ref{s-indiv}, we write $q_i=p^1_i$ and in order for our problem
to be feasible, we assume $0\leq p < q \leq p+K$.

Let us show that the LP in (\ref{eq:lp-tight}) is feasible.
Indeed, we can form a piecewise affine function that is feasible
for (\ref{eq:infinite-LP}) by taking $C(w_{0},K_{0})=E_{\pi
}(w_{0}^{T}x-K_{0})_{+}$, where $\pi $ is the probability measure
defined in (\ref{eq:dist-feas}). By construction, this function
corresponds to a feasible point of (\ref{eq:lp-tight}) and the
variables $g_{i}$ in (\ref{eq:lp-tight}) are simply the
subgradients of $C(w_{0},K_{0})$ at the data points. Finally, the
LP in (\ref{eq:lp-tight}) is finite, since we always have $0\leq
E_{\pi }(w_{0}^{T}x-K_{0})_{+}\leq w_{0}^{T}q$ and the feasible
set of (\ref{eq:lp-tight}) is compact. This means that the optimum
in (\ref{eq:lp-tight}) is attained.

Suppose that the forward price information is ignored, i.e.
assuming that $m=n$, and $w_{0}\in \mathbf{R}_{+}^{n}$. We note
$e_{i}$, the $i$-th unit vector. Without loss of generality, we
set $w_{0}^{T}e=1$. Since the function
$C(w_{0},K_{0})=w_{0}^{T}p+\left( w_{0}^{T}K-K_{0}\right) _{+}$ is
a feasible point of the infinite LP in (\ref{eq:infinite-LP}), if
we call $V^{\mathrm{LP}}$ the upper bound computed by the linear
program (\ref {eq:lp-tight}), we must have:
\[
V^{\mathrm{LP}}\geq w_{0}^{T}p+\left( w_{0}^{T}K-K_{0}\right) _{+}.
\]
Now, using the necessary conditions in (\ref{eq:infinite-LP}) and the
convexity of
\[
{\mathbf{E}}_{\pi _{\varepsilon }}\left( w^{T}x-K\right) _{+}
\]
in $(w,K)$ we can write
\begin{eqnarray*}
{\mathbf{E}}_{\pi _{\varepsilon }}\left( w_{0}^{T}x-K_{0}\right)
_{+} &=&{\mathbf{E}}_{\pi _{\varepsilon }}\left( w_{0}^{T}x-\left(
w_{0}^{T}K+\left(
K_{0}-w_{0}^{T}K\right) \right) \right) _{+} \\
&\leq &\sum_{i=1}^{n}w_{0,i}{\mathbf{E}}_{\pi _{\varepsilon }}\left(
x_{i}-\left( K_{i}+\left( K_{0}-w_{0}^{T}K\right) \right) \right) _{+} \\
&=&\sum_{i=1}^{n}w_{0,i}{C}\left( e_{i},K_{i}+\left( K_{0}-w_{0}^{T}K\right)
\right) .
\end{eqnarray*}
The conditions on the slope of the function $C(w,K)$ imply
\[
\sum_{i=1}^{n}w_{0,i}{C}\left( e_{i},K_{i}+\left( K_{0}-w_{0}^{T}K\right)
\right) \leq w_{0}^{T}p+\left( w_{0}^{T}K-K_{0}\right) _{+}.
\]
Hence, $V^{\mathrm{LP}}\leq w_{0}^{T}p+\left( w_{0}^{T}K-K_{0}\right) _{+}$
and finally
\begin{equation}
V^{\mathrm{LP}}=w_{0}^{T}p+\left( w_{0}^{T}K-K_{0}\right) _{+},
\label{eq:LP-sharpness}
\end{equation}
where we recover the expression found in (\ref{eq:upper-no-forward-bnd}).
This means that the upper bound computed by the LP relaxation is tight in
the particular case considered above.

Now we turn to the case when forward price constraints
$\mathbf{E}_{\pi }x_{i}=q_{i}\ $ for $i=1,\ldots ,n,$ are
included. As already observed in \ref{s-upper-forwards}, the
function
\[
d^{\mathrm{sup}}(w_{0},K_{0})=\displaystyle\max_{0\leq j\leq
n+1}:w_{0}^{T}p+\sum_{i}w_{0,i}\displaystyle\min (q_{i}-p_{i},\beta
_{j}K_{i})-\beta _{j}K_{0},
\]
is convex in $(w_{0},K_{0}).$ Also, when $w_{0}=e_{i}$, and
$K_{0}=K_{i}$, we obtain $d^{\mathrm{sup}}=p_{i}$, while for
$K_{i}=0$, we obtain $d^{\mathrm{sup}}=q_{i}$. This means that
$d^{\mathrm{sup}}(w,K)$ is a feasible point of the infinite
program (\ref{eq:infinite-LP}) and hence $V^{\mathrm{LP}}\geq
d^{\mathrm{sup}}(w_{0},K_{0})$.

Since the finite LP (\ref{eq:lp-tight}) is attained, at a point
denoted by
\[
z^{\ast }=\left[ p_{0}^{\ast },g_{0}^{\ast T},\ldots ,g_{k}^{\ast T}\right]
^{T},
\]
we can define the call price function
\[
d^{\mathrm{LP}}(w,K)=\max_{i=0,\ldots ,2n+1}\left\{ p_{i}
+\left\langle g_{i}^{\ast },(w,K)-(w_{i},K_{i})\right\rangle
\right\} ,
\]
with $p_{0}=p_0^{\ast}$. We can compute the corresponding strike
prices $\hat{K}=t^{\ast }K$ and option prices $\hat{p}=(1-t^{\ast
})q+t^{\ast }p$, as in \ref{s-upper-forwards-proba}. By convexity
of $d^{\mathrm{LP}}(w,K),$ we have
$d^{\mathrm{LP}}(e_{i},\hat{K})\leq \hat{p}_{i}$ for $i=1,\ldots
,n$. From (\ref{eq:LP-sharpness}), we get that
$d^{\mathrm{LP}}(w_{0},K_{0})=V^{\mathrm{LP}}\leq
d^{\mathrm{sup}}(w_{0},K_{0})$, hence finally
$d^{\mathrm{LP}}(w_{0},K_{0})=d^{\mathrm{sup}}(w_{0},K_{0})$. This
shows that the LP relaxation of the upper bound is tight when the
input is composed of forwards and one option price per asset.

Now, from section \ref{ss:two-options}, we know that the upper
bound problem given two option prices per asset can be reduced to
an upper bound problem given forwards and one option price per
asset. This means that the LP relaxation of the upper bound
problem is also tight in the case where two option prices are
given for each asset. Finally, we conclude by remarking that
\cite[Th. 4.1]{Hobs04} show that the optimal upper bound involves
at most two option prices per asset, hence the problem of finding
an upper bound in the general problem (\ref{eq:equity-prob}),
\ie~given option prices for many strikes, reduces to an upper
bound problem given only two options per asset, for which the LP
relaxation is tight.
\end{proof}

The relaxation we obtain in section 2 treats the upper and lower
bound problems in a completely symmetric way. This does not
\emph{at all} reflect the inherent complexity of the exact
problems. In the upper bound case, we are essentially lucky and we
can find, in most cases, discrete distributions with compact
support matching the upper bounds, hence we very often get perfect
duality (most notably in the special case discussed in this
section). This is not true for the lower bound and we almost never
get perfect duality in general.

We can also remark that since the largest decreasing convex
functions interpolating the market option prices are piecewise
affine, the upper bound computed in \cite[Theorem 4.1]{Hobs04} can
also be formulated as the solution to a linear program.

\section{Numerical results}
\label{s-numerical} In this section, we first show how a variant
to the program in Proposition \ref{prop:relax-general} can
efficiently ``clean'' the market data on single asset option
prices to remove the non convexities due to noise (asynchronous
data, liquidity, etc). We then focus on an equity market example
and compare the upper bound computed using the linear programming
relaxation in Proposition \ref{prop:tightness} and the bounds
obtained in \cite{Hobs04} on the data set described in \cite[Table
2]{Hobs04}. We show that this same relaxation is not tight in the
lower bound case by exhibiting an example where the lower bound
computed in section \ref{ss:lower-bound} is larger than the
relaxation's result. Finally, we test our relaxation technique in
the general case using a simplified interest rate data set taken
from \cite{Brig01} and show how the lower and upper bounds perform
in the general setting of computing bounds on the price of a
basket given other basket prices.

\subsection{Basket bounds given single asset option prices}
\subsubsection{Market data}
\begin{figure}[ht]
\begin{center}
\psfrag{str}[t][b]{Strike Price}
\psfrag{pr}[b][t]{Option Price}
\psfrag{DOW}[c][t]{Basket Option Prices}
\includegraphics[width=0.7 \textwidth]{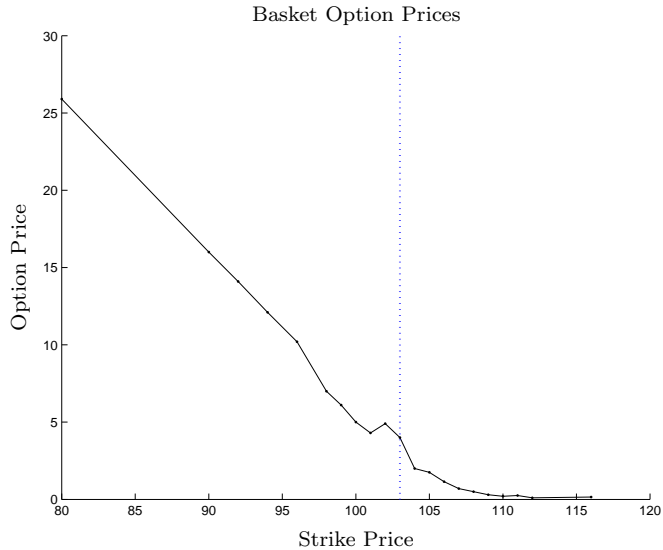}
\caption{\label{fig:dowprices} Dow Jones index call option prices
on Mar. $17$ $2004$, with maturity Apr. $16$ $2004$, for various
strikes. The vertical dotted line is at the forward price.}
\end{center}
\end{figure}

In Figure \ref{fig:dowprices}, we show the market price on Mar.
$17$, $2004$ of a certain number of basket call options on the Dow
Jones index, with maturity Apr. $16$, $2004$, for various strikes
(the underlying asset of this option is here the Dow Jones index
divided by $100$).
We notice that the convexity requirement with respect to the
strike price in Proposition~\ref{prop:nec-suff-convex} is violated
between strikes $95$ and $103$ for example. Since prices in this
plot are last quotes, this does not necessarily mean that an
arbitrage is present. In practice, options with strikes far away
from the forward price (here equal to $103$) tend to be somewhat
illiquid and their last quoted price can be an unreliable
indication of the price at which they would trade today. However,
since this option price data is used to calibrated pricing models,
being able to efficiently detect such anomalies in the data is
crucial. We can use the result of
Proposition~\ref{prop:relax-general} to clean the input data on
single asset option prices, from \cite{Laur00} we know that the
relaxation in Proposition \ref{prop:relax-general} is always exact
in dimension one and we look for a decreasing, convex
approximation of the market data. Let $p_i$ be the market prices
of calls on an asset $x$ with strikes $K_i$, for $i=1,\ldots,m$,
we can find the closest (in the $l_1$ norm sense) arbitrage free
approximation of the market prices by solving:
\begin{equation}
\label{eq:clean-data}
\begin{array}{ll}
\mbox{minimize} & \sum_{i=1}^m |y_i-p_i|\\
\mbox{subject to} & \langle g_i,(0,K_j-K_i)\rangle \leq y_j-y_i\\
& g_{i,1} \geq 0,~-1 \leq g_{i,2} \leq 0\\
& \langle g_i,(1,K_i)\rangle=y_i,\quad i,j=1,\ldots,m,
\end{array}
\end{equation}
which is a linear program in the variables $g_i \in \reals^2$ and
$y_i\in\reals$, for $i=1,\ldots,m$.

\subsubsection{Upper bound with many strikes}
Here we test empirically the result detailed in Proposition
\ref{prop:tightness} using the data set in \cite[Table 2]{Hobs04}.
This data set still contains minor non convexities (on the BA
option prices for example) and we first clean it using the
procedure described in (\ref{eq:clean-data}) above. We then
compare the upper bound on the price of a basket option with
uniform weights 0.071 obtained using the relaxation in Proposition
\ref{prop:relax-general} with those detailed in \cite[Table
3]{Hobs04}. The results are detailed in Table
\ref{tab:upper-hobs04}. As expected from the result in Proposition
\ref{prop:tightness}, the bounds match up to a small numerical
error that is most likely due to differences in the data cleaning
procedure.

\begin{table}
{\small\[
\begin{array}{cccc}
\mbox{\bf ~ Strikes ~} & \mbox{\bf ~DJX price~} & \mbox{\bf ~LP
relax.~} & \mbox{\bf ~U.B.~}\\
\hline\\
52  &  47.10 & 47.13 & 47.09\\
56  & 43.10 & 43.14 & 43.10\\
60  & 39.10 & 39.15 & 39.11\\
64  & 35.10 & 35.16 & 35.11\\
68  & 31.10 & 31.17 & 31.12\\
70  & 29.10 & 29.18 & 29.13\\
72  & 27.10 & 27.19 & 27.14\\
76  & 23.10 & 23.21 & 23.15\\
80  & 19.10 & 19.25 & 19.18\\
84  & 15.20 & 15.37 & 15.24\\
88  & 11.30 & 11.63 & 11.42\\
90  & 9.40 & 9.80  & 9.61\\
92  & 7.50 & 8.08  & 7.90\\
94  & 5.80 & 6.49  & 6.32\\
95  & 4.95 & 5.73  & 5.57\\
96  & 4.15 & 5.01  & 4.85\\
97  & 3.35 & 4.33  & 4.19\\
98  & 2.73 & 3.71 & 3.58\\
99  & 2.13 & 3.14 & 3.02\\
100  & 1.60 & 2.64 & 2.53\\
102  & 0.78 & 1.82 & 1.73\\
103  & 0.50 & 1.48 & 1.42\\
104  & 0.33 & 1.22 & 1.16\\
105  & 0.15 & 1.00 & 0.95\\
106  & 0.15 & 0.80 & 0.75\\
107  & 0.15 & 0.64 & 0.59\\
\end{array}
\]}
\caption{\label{tab:upper-hobs04} This table displays for various
strikes: the market price (DJX) of a basket option with uniform
weights, the upper bound (LP relax.) computed using the relaxation
in $\S\ref{s-general}$ and the upper bound (U.B.) computed in
\cite{Hobs04}.}
\end{table}

Intuitively, tightness in this setting stems from the fact that
the convex function $C(w,K)$ which is the optimal solution to the
relaxation in Proposition \ref{prop:relax-general} coincides on
each axis with the convex envelope of the single asset option
prices in the data set (the functions $\bar C^{(i)}$ in \cite[p.
12]{Hobs04}). Unfortunately, this does not happen in the more
general case where the data consists of basket option prices.

\subsubsection{Lower bound}
We can also test the tightness of the various bounds obtained
above on a simulated data set and compare these bounds with actual
prices. The forward prices are given by $q=\{70, 50, 40, 40,
40\}$. We set $K=\{70,50,40,40,40\}$ and
$p=\{1.61,1.43,0.93,0.70,0.47\}$. Using the results of section
\ref{s-indiv}, we get bounds on the price of a basket option with
weight vector $w_{0}=\{0.2,0.2,0.2,0.2,0.2\}$ which are detailed
in Table \ref{tab:lower}.
\begin{table}
{\small\[
\begin{array}{rccccc}
\mbox{\bf Strike price} & 3.84 & 4.32 & 4.80 & 5.28 & 5.76\\
\hline\\
\mbox{\bf Upper bound (relax.)} & 1.71 & 1.37 & 1.03 & 1.03 & 1.03\\
\mbox{\bf Upper bound (from $\S\ref{s-indiv}$)} & 1.71 & 1.37 & 1.03 & 1.03 & 1.03\\
\mbox{\bf Lower bound (from $\S\ref{s-indiv}$)} & 0.96 & 0.48 & 0.09 & 0.00 & 0.00\\
\mbox{\bf Lower bound (relax.)} & 0.96 & 0.48 & 0.00 & 0.00 & 0.00\\
\end{array}
\]}
\caption{\label{tab:lower} Lower and upper bounds on the price of
a basket given single asset forward prices and one option price
per asset. We compare the lower bounds produced in
$\S\ref{s-indiv}$ with those produced by the linear programming
relaxation in $\S\ref{s-general}$.}
\end{table}

We notice that the lower bound computed using
(\ref{eq:LP-with-forwards}) is tighter than that provided by the
LP relaxation in (\ref{eq:finite-LP}) for at least one strike
price. In this case, the LP relaxation is equal to the trivial
lower bound given by the intrinsic value of the option.

\subsection{Basket bounds given basket option prices}
Here, we test the relaxation in Proposition
\ref{prop:relax-general} on a simplified data set taken from
\cite{Brig01}. Our objective here is to look for bounds on the
price of a swaption given other swaption prices in a simplified
setting. We refer the reader to \cite{Rebo98} and \cite{dasp02d}
for further details.

We compute market option prices by simulation. The underlying
asset prices follow lognormal dynamics as in \cite{Blac73}. The
forwards are given here by $F=\{0.03,0.04,0.04,0.05,0.05\}$ (these
are forward interest rates), and the covariance matrix C is taken
from \cite[pp. 301,~311]{Brig01}: {\small
\[ C=\left(
\begin{array}{lllll}
0.034 & 0.032 & 0.026 & 0.021 & 0.018 \\
0.032 & 0.035 & 0.019 & 0.026 & 0.011 \\
0.026 & 0.019 & 0.024 & 0.010 & 0.019 \\
0.021 & 0.026 & 0.010 & 0.020 & 0.004 \\
0.018 & 0.011 & 0.019 & 0.004 & 0.017
\end{array}
\right)
\]}

\begin{figure}[ht]
\begin{center}
\psfrag{k}[t][b]{Strike Price}
\psfrag{vol}[b][t]{Implied Volatility}
\includegraphics[width=0.9 \textwidth]{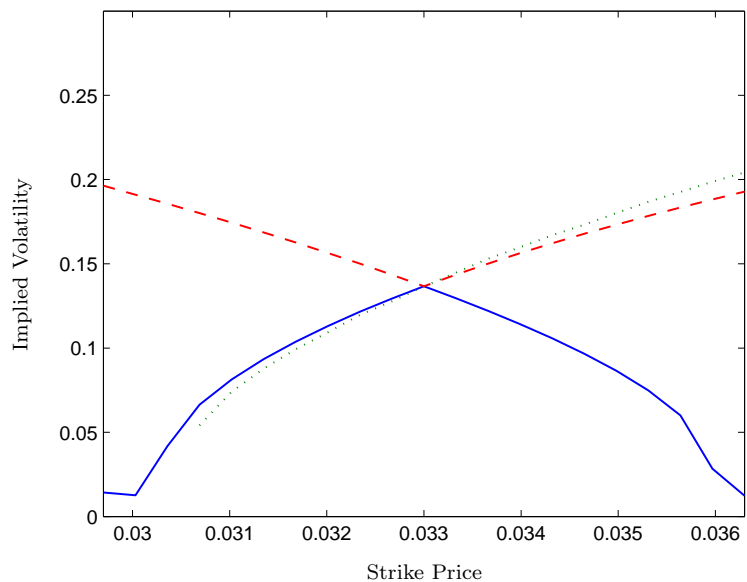}
\end{center}
\caption{\label{fig:figbounds} Upper (dashed) and lower (solid)
bounds on the price (plotted here in terms of implied volatility)
of a basket option given other basket option prices. These bounds
are computed using the linear programming relaxation detailed in
$\S\ref{s-general}$. The price of the at the money option is given
here, hence the upper and lower bounds match at this point. The
dotted line is the market volatility.}
\end{figure}

We look for upper and lower bounds on the price of a basket option
(a yield curve option here) with weight vector
$w_{0}=\{.2,.1,.2,.1,.2\}$ for various strikes, given the
simulated price of other at the money basket options (swaptions)
with weights vectors given by the lines of the matrix: {\small\[
W=\left(
\begin{array}{lllll}
1 & 0 & 0 & 0 & 0\\
1 & 1 & 0 & 0 & 0\\
1 & 1 & 1 & 0 & 0\\
1 & 1 & 1 & 1 & 0\\
1 & 1 & 1 & 1 & 1\\
\end{array} \right)
\]}
and the at the money price of a basket option with weights $w_0$
(which is usually liquid). We are looking for bounds on basket
options with weight $w_0$ and strike prices $\pm 10\%$ away from
the forward. The resulting price bounds are plotted in Figure
\ref{fig:figbounds} in terms of implied volatility. We notice that
the market volatility is actually equal (up to a small numerical
error) to either the lower bound when the option is in the money
or to the upper bound when the option is out of the money.

\bibliographystyle{amsalpha}
\bibliography{MainPerso}
\end{document}